\documentclass[12pt]{article}%
\usepackage{amsmath}
\usepackage{amsfonts}
\usepackage{amssymb}
\usepackage{graphicx}%
\providecommand{\U}[1]{\protect\rule{.1in}{.1in}}
\newtheorem{theorem}{Theorem}

\newtheorem{corollary}[theorem]{Corollary}

\newtheorem{definition}[theorem]{Definition}

\newtheorem{lemma}[theorem]{Lemma}

\newtheorem{proposition}[theorem]{Proposition}
\newtheorem{remark}[theorem]{Remark}

\newenvironment{proof}[1][Proof]{\noindent\textbf{#1.} }{\ \rule{0.5em}{0.5em}}
\begin{document}

\title{The sup-completion of a Dedekind complete vector lattice.}
\author{Youssef Azouzi and Youssef Nasri\\{\small Research Laboratory of Algebra, Topology, Arithmetic, and Order}\\{\small Department of Mathematics}\\{\small Faculty of Mathematical, Physical and Natural Sciences of Tunis}\\{\small Tunis-El Manar University, 2092-El Manar, Tunisia}}
\date{}
\maketitle

\begin{abstract}
Every Dedekind complete Riesz space $X$ has a unique sup-completion $X^{s}$,
which is a Dedekind complete lattice cone. This paper aims to present a
systematic study this cone by extending several known results to general
setting, proving new results and, in particular, introducing for elements of
$X^{s}$ finite and infinite parts. This enuables us to get a satisfactory
abstract formulation of some classical results in the setting of Riesz spaces.
We prove, in pareticular, a Riesz space version of Borel-Cantelli Lemma and
present some applications to it.

\end{abstract}

\section{Introduction}

Completeness is a desirable property for anyone dealing with Riesz spaces.
Such spaces may have several kinds of completeness. Under laterally
completeness suprema exist for all disjoint families. Under Dedekind
completeness suprema exist for all nonempty bounded subsets from above.
Dedekind complete Riesz spaces are also called order complete; this can be
explained by the fact that in those spaces order Cauchy nets are order
convergent. If we consider uniform convergence instead of order convergence we
talk about uniformly completeness. Recently, the notion of unbounded order
convergence (uo-convergence) has received much attention (see for example
\cite{L-174,L-65,L-444,L-745}) and then one can be interested in uo-complete
Riesz spaces. Uo-completeness means that every uo-Cauchy net is uo-convergent.
It is quite surprising that Uo-completeness is equivalent to universal
completeness. This is the main result in \cite[Theorem 17]{L-444}. When nets
are replaced by sequences we speak about sequential completeness. In general,
sequential completeness is not equivalent to completeness. If the spaces we
consider fail to be complete there is a way to make them complete by embedding
those spaces in complete spaces (of course, under some additional
assumptions). The notion of completion is then helpful and useful. Among
several types of completion in Riesz space Theory we will focus in this paper
on the notion of sup-completion. This notion has been introduced by Donner in
\cite{b-1665}; it is perhaps less known and less frequently used in the
literature. Several recent papers present important applications of it,
essentially published by Grobler \cite{L-06,L-745}, Grobler and Labushangne
\cite{L-12,L-13}, and the first author \cite{L-444}. Sup-completion is used by
Grobler in \cite{L-06} to constructed Daniell integral in Riesz spaces.
Recently the first author used this notion to prove the equivalence between
universal completion and unbounded order completion. The current paper can be
broken into two parts, both of them deal with sup-completion. The first part
is devoted to the study of the sup-completion of a Dedekind complete Riesz
space itself. It starts by a brief review of what is known and then presents
several new results, namely the introduction of finite and infinite parts of
elements of $X^{s}.$ The second part is motivated by the famous Borel-Cantelli
Lemma. It aims to give Riesz space generalizations of this lemma and provide
some applications. The Borel-Cantelli Lemma contains two parts. Its first one
is stated under the assumption that the series $%
{\textstyle\sum\limits_{k}}
\mathbb{P}\left(  A_{k}\right)  $ is convergent, where $\left(  A_{n}\right)
$ is a sequence of events, and its generalization to the setting of Riesz
spaces is quite obvious. This has already been done by Kuo, Labuschagne and
Watson in \cite{L-260}. Another generalization of the first Borel-Cantelli
lemma, due to Barndorff-Nielsen and Balakrishnan-Stepanov, have been recently
translated to the setting of Riesz spaces by Mushambi, Watson and Zinsou
\cite{L-638}. The second part of the Borel Cantelli Lemma, however, has not
received enough attention in the frame of Riesz spaces setting and it is
obvious that its generalization requires more care. The condition $%
{\textstyle\sum\limits_{k}}
\mathbb{P}\left(  A_{k}\right)  =\infty$ forces us to work in the space
$X^{s}.$ A good understanding of that cone is then required in order to get
satisfactory translation of several results from classical theory of
probability to the setting of Riesz spaces. A series in a Riesz spaces may be
converging on some band $B$ and diverging on its orthogonal $B^{d}$ and we can
be interested in determining the largest band on which the series is
converging. So, we are led to consider the finite and infinite parts of an
element in $X^{s}.$ This allows us to get more meaningful statements in the
case of Riesz spaces.

An outline of the paper is as follows: Section 2 contains some preliminaries
on several types of convergence considered in Riesz spaces. Section 3 is
devoted to the notion of sup-completion: we will describe the basic properties
of the cone $X^{s}$ and extend several of them. We introduce for an element
$x$ in $X_{+}^{s}$ its finite and infinite parts and use these notions to
extend Borel-Cantelli Lemma in Section 4. In this section we present a new
version of the first Borel Cantelli Lemma (BCL1), which extends the one
obtained earlier by Kuo, Labushagne and Watson and prove a Riesz space version
of the second Borel-Cantelli Lemma (BCL2). The last section, Section 5,
provides several applications to the previous sections and especially to
Borel-Cantelli Lemma.

\section{Preliminaries}

Throughout $X$ denotes a Dedekind complete Riesz space. So, every band in $X$
is a a projection band. For a band $B$ we denote by $P_{B}$ the associated
band projection. If $B=B_{x}$ is a principal band generated by $x,$ we write
$P_{x}$ instead of $P_{B_{x}}.$ If $P$ is a band projection we let
$P^{d}:=I-P$ denote the band projection on the band $B^{d}$. The universal
completion of $X$ is denoted by $X^{u},$ while its sup-completion is denoted
by $X^{s}$.

We will deal in this paper with three modes of convergence; each defines a
linear convergence structure on $X$ in the sense of \label{L-748}. We briefly
recall their definitions and some basic and important facts about them that
will be useful for us. The first one, and the most useful, is the order convergence.

\begin{definition}
We say that a net $\left(  x_{\alpha}\right)  _{\alpha\in A}$ in a vector
lattice $X$ is order convergent to $x$ if there exists a `dominating' net
$\left(  y_{\beta}\right)  _{\beta\in B}$ satisfying $y_{\beta}\downarrow0$
and for any $\beta\in B$ there exists $\alpha_{0}\in A$ such that $\left\vert
x_{\alpha}-x\right\vert \leq y_{\beta}$ for every $\alpha\geq\alpha_{0}.$ We
write $x_{\alpha}\overset{o}{\longrightarrow}x$ or, more simply, $x_{\alpha
}\longrightarrow x.$
\end{definition}

The second mode of convergence can be viewed as an abstraction of almost
surely convergence.

\begin{definition}
We say that a net $\left(  x_{\alpha}\right)  $ in a Riesz space $X$ unbounded
order converges (or, uo-converges) to $x,$ and we write $x_{\alpha}%
\overset{uo}{\longrightarrow}x,$ if for every $u\in X_{+},$ the net
$\left\vert x_{\alpha}-x\right\vert \wedge u$ is order convergent to $0.$
\end{definition}

If $X$ has a weak order unit $e$ then $x_{\alpha}\overset{uo}{\longrightarrow
}x$ if and only if $\left\vert x_{\alpha}-x\right\vert \wedge e\overset
{o}{\longrightarrow}0.$

It is clear that order convergence agrees with uo-convergence for eventually
bounded nets. They agree also for sequences if the space is universally
complete but they do not for nets. We mention also a very useful result
obtained in \cite{L-65}: If $X$ is a vector lattice, $Y$ a regular vector
sublattice of $X$ and $\left(  y_{\alpha}\right)  $ a net in $Y$ then
$y_{\alpha}\overset{uo}{\longrightarrow}0$ in $Y$ if and only if $y_{\alpha
}\overset{uo}{\longrightarrow}0$ in $X.$ This can be applied, in particular,
when $X=Y^{u}$ is the universal completion of a Dedekind complete Riesz space
$Y$ since in this case $Y$ is an ideal of $X$ and every ideal is regular.

The third mode of convergence can be viewed as a generalization of convergence
in probability. We assume here that $X$ is equipped with a conditional
expectation operator $T$ with $Te=e.$ Recall that $T$ is an order continuous
strictly positive projection which has $R\left(  T\right)  $ a Dedekind
complete Riesz subspace.

\begin{definition}
Let $E$ be a Dedekind complete Riesz space with weak order unit $e,$ and let
$T$ be a conditional expectation operator on $E$ satisfying $Te=e.$ We say
that a net $(x_{\alpha})_{\alpha\in A}$ converges in $T$-conditional
probability to $x,$ and we write $x_{\alpha}\overset{TP}{\longrightarrow}x,$
if the net $\left(  TP_{(|x_{\alpha}-x|-\epsilon e)^{+}}e\right)  $ converges
in order to $0$ for each $\epsilon>0$.
\end{definition}

This notion has been introduced in \cite{L-14} as a generalization of
convergence in probability. Indeed, it agrees with the convergence in
probability for sequences in the case when $X=L_{1}\left(  \Omega
,\mathcal{F},\mathbb{P}\right)  $ and $T=\mathbb{E}$ is the expectation
operator. It is perhaps worth mentioning that if $T=\operatorname*{Id}%
\nolimits_{X}$ is the identity map, then the convergence in $T$-conditionally
probability is the unbounded order convergence. So, the following result is a
generalization of \cite[Corollary 3.5]{L-65}; it could also be compared with
\cite[Lemma 2.11]{L-173}.

\begin{lemma}
\label{LA}Let $X$ be a Dedekind complete Riesz space with weak order unit $e$
and $T$ a conditional expectation with $Te=e.$ For a net $(x_{\alpha}%
)_{\alpha\in A}$ in $X$, the following are equivalent.

\begin{enumerate}
\item[(i)] $x_{\alpha}\longrightarrow x$ in $T$-conditional probability;

\item[(ii)] $T(|x_{\alpha}-x|\wedge u)\overset{o}{\longrightarrow}0$ for every
$u\in X_{+};$

\item[(iii)] $T(|x_{\alpha}-x|\wedge e)\overset{o}{\longrightarrow}0.$
\end{enumerate}
\end{lemma}

\begin{proof}
Without loss of generality we may assume that $x=0.$

(i) $\Longrightarrow$ (iii) Let $\varepsilon\in\left(  0,\infty\right)  .$and
put $P_{\varepsilon}=P_{\left(  |x_{\alpha}|\wedge e-\varepsilon e\right)
^{+}}.$ Since $P_{\varepsilon}^{d}\left\vert x_{\alpha}\right\vert
\leq\varepsilon e$ we have%
\begin{align*}
|x_{\alpha}|\wedge e  &  =P_{\varepsilon}(|x_{\alpha}|\wedge e)+P_{\varepsilon
}^{d}(|x_{\alpha}|\wedge e)\\
&  \leq P_{\varepsilon}e+P_{\varepsilon}^{d}(|x_{\alpha}|)\leq P_{\varepsilon
}e+\varepsilon e.
\end{align*}
Now apply $T$ to the above display and taking the limit supremum over $\alpha$
to get%
\[
\limsup\limits_{\alpha}T(|x_{\alpha}|\wedge e)\leq T\varepsilon e=\varepsilon
e.
\]
As this happens for every $\varepsilon>0$ we derive that $\limsup
\limits_{\alpha}T(|x_{\alpha}|\wedge e)=0$ and then $\lim T(|x_{\alpha}|\wedge
e)=0$ as required.

(iii) $\Longrightarrow$ (i) Assume that $T(|x_{\alpha}|\wedge e)\overset
{o}{\longrightarrow}0.$ Then for every $\varepsilon\in\left(  0,1\right)  $ we
have $P_{\left(  |x_{\alpha}|-\varepsilon e\right)  ^{+}}e\leq e\leq
\varepsilon^{-1}e$. On the other hand, from the inequality%
\[
P_{\left(  |x_{\alpha}|-\varepsilon e\right)  ^{+}}e\leq\varepsilon^{-1}|x|,
\]
it follows that%
\[
TP_{\left(  |x_{\alpha}|-\varepsilon e\right)  ^{+}}e\leq\varepsilon
^{-1}T(|x_{\alpha}|\wedge e)\overset{o}{\longrightarrow}0,
\]
\newline which gives (ii).

(ii) $\Longleftrightarrow$ (iii). The forward implication is trivial. For the
converse assume that (iii) occurs and let $u\in X_{+}.$ Then%
\begin{align*}
T(|x_{\alpha}|\wedge u)  &  =T(|x_{\alpha}|\wedge u-(|x_{\alpha}|\wedge
u\wedge ke)+T(|x_{\alpha}|\wedge e)\\
&  \leq T(u-u\wedge ke)+T(|x_{\alpha}|\wedge ke).
\end{align*}
Hence $\limsup\limits_{\alpha}T(|x_{\alpha}|\wedge u)\leq T(u-u\wedge ke)$ for
every $k>0.$ Letting $k$ to $\infty$ yields $\limsup\limits_{\alpha
}T(|x_{\alpha}|\wedge u)=0,$ which shows that $T(|x_{\alpha}|\wedge
u)\overset{o}{\longrightarrow}0$ as required.
\end{proof}

For each mode of conververgence mentioned above we may define Cauchy net to be
a net $\left(  x_{\alpha}\right)  _{\alpha\in A}$ such that the net $\left(
x_{\alpha}-x_{\beta}\right)  _{\left(  \alpha,\beta\right)  \in A\times A}$
converges to $0$.

For more information about order convergence and unbounded order convergence
the reader is referred to \cite{L-65,L-444} and references therein and to
\cite{L-14} for convergence in $T$-conditionally probability.

\section{Properties of the sup-completion}

We consider in this section a Dedekind complete Riesz space $X.$ The notion of
sup-completion has been introduced by Donner in \cite{b-1665} and it is used
there to prove some extension Theorems. It has been explored by Grobler in
\cite{L-06} to construct Daniell integral and then develop a kind of
Functional Calculus in Riesz spaces. Some extra results have been obtained by
the first author in \cite{L-444} and used to get a presentation theorem of
element of $X^{s}$ by integrals. This is a crucial step to prove the main
theorem in \cite{L-444}, which states that a Riesz space is universally
complete if and only if it is uo-complete. To begin the discussion, let us
recall briefly the construction of the sup-completion. Consider the set
$\mathcal{A}$ of all nonempty, upward directed subsets of $X$ endowed with the
equivalence relation $\sim$ given by%
\[
A\sim B\Longleftrightarrow\sup\limits_{a\in A}\left(  x\wedge a\right)
=\sup\limits_{b\in B}\left(  x\wedge b\right)  \text{ for all }x\in X.
\]
Addition, multiplication by positive scalars, and ordering are defined on
$\mathcal{A}$ in a natural way. Endowed with the quotient operations, the
quotient space $\mathcal{A}\diagup\sim$\ is an ordered cone, called the
sup-completion of $X$ and denoted by $X^{s}.$ The space $X$ can be viewed as a
subset of $X^{s}$ by identifying $x$ with the class of $\left\{  x\right\}  $
for each element $x\in X.$Theorem 1.4 in \cite{b-1665} provides a
characterization of the cone $X^{s}$ and contains the fundamental properties
of that cone. To make the reading of this paper easy we list below the most
important properties of $X^{s},$ which include the extra results obtained in
\cite{L-444}. The reader can observe that some of these properties are stated
with different assumptions, which is not easy to keep in mind (see for
example, (P3) and (P6) below) and it is natural to ask whether these
properties hold without any extra assumption. Our first purpose in this
section is to extend these properties in more general setting by relaxing the
assumptions, which makes them more natural and makes their utilization easier.
The second purpose is to introduce finite and infinite parts of elements of
$X^{s}$ and investigate their properties. This enables us to get satisfactory
abstract formulation of some results in the setting of Riesz spaces. It is the
case of the second Borel-Cantelli lemma (Theorem \ref{BCL2}) and Theorems
\ref{T2} and \ref{T3}. For any element $y\in X^{s},$ let $\left[  y\right]
^{\leq}$ denote the subset $\left\{  x\in X:x\leq y\right\}  .$ We are now
ready to list several important properties of the cone $X^{s}.$ The reader is
referred to \cite{b-1665} and \cite{L-444} for their proofs.

\begin{enumerate}
\item[(P1)] $X$ is the set of invertible elements in $X^{s}$ with coinciding
algebraic and order structures.

\item[(P2)] For every $y\in X^{s}$ we have $y=\sup\left[  y\right]  ^{\leq}$.

\item[(P3)] For every $x,\ y\in X^{s},\ a\in X$ we have
\[
x+a\wedge y=(x+a)\wedge(x+y).
\]

\item[(P4)] If $x\in X$ and $y\in X^{s}$ satisfy $y\leq x$ then $y\in X$.

\item[(P5)] $X^{s}$ has a greatest element.

\item[(P6)] For any two non-empty subsets $A,\ B\subset X^{s}$ satisfying
$\sup A=\sup B$ the equality%
\[
\sup_{a\in A}(a\ \wedge x)=\sup_{b\in B}(b\wedge x)
\]
holds for every $x\in X$. In particular, if $m=\sup A$ then%
\[
m\wedge x=\sup_{a\in A}(a\ \wedge x).
\]

\item[(P7)] If $X$ has a weak order unit $e$ and $0\leq x\in X^{s}$, then
\[
x=\sup\limits_{k\geq1}(ke\wedge x).
\]

\item[(P8)] If $A\ $and $B$ are non-empty subsets of $X$ then
\[
\sup(A+B)=\sup A+\sup B.
\]

\item[(P9)] \textit{Birkhoff Inequality}: If $a,b\in X,$ $c\in X^{s}$ then%
\[
|a\wedge c-b\wedge c|\leq|a-b|.
\]

\item[(P10)] \textit{Riesz Decomposition Property:} \textit{If} $0\leq
x,\ y,\ z\in X^{s}$ \textit{with} $x\leq y+z$ \textit{then there exist}
$y_{1},\ z_{1}\in X^{s}$ \textit{such that} $y_{1}\leq y,\ z_{1}\leq z$
\textit{and \qquad}$x=y_{1}+z_{1}$.
\end{enumerate}

Properties (P1)-(P8) were proved by Donner in \cite{b-1665}, while properties
(P9), (P10) have been recently shown by the first author in \cite{L-444}. The
greatest element of $X$ will be denoted by $\infty.$

As we have already mentioned above, the sup-completion may be a suitable space
when one wants to extend maps, especially those extensions that preserve
supremum for increasing nets. This occurs when dealing with integrals of
functions taking their values in $\left[  0,\infty\right]  $ and allow
integrals to be infinite. The following result goes in this direction. It
extends \cite[Proposition 3]{L-444} and improves it.

\begin{theorem}
\label{P1}\textit{Let} $X$ \textit{and} $Y$ \textit{be two Dedekind complete
Riesz spaces and }$f:X\rightarrow Y$ \textit{be an order continuous increasing
map. Then} the following statements hold.

\begin{enumerate}
\item[(i)] $f$ \textit{can be extended to a unique left order continuous
increasing map} $f^{s}$ \textit{from} $X^{s}$ \textit{to} $Y^{s}$. \textit{If}
$f$ \textit{is additive} (\textit{resp. positively homogeneous}),
\textit{then} so is $f^{s}$.

\item[(ii)] \textit{If} $f$ is only defined from $X_{+}$ to $Y_{+}$ then $f$
has a unique \textit{left order continuous extension }$f^{s}$ \textit{from}
$X_{+}^{s}$ \textit{to} $Y_{+}^{s}$ which is increasing. \textit{Moreover, if}
$f$ \textit{is additive} (\textit{resp. positively homogeneous}), \textit{then
so is} $f^{s}.$

\item[(iii)] If $f$ is, in addition, a linear projection and $f\left(
X\right)  $ is a Riesz subspace of $Y=X,$ then $f\left(  X\right)  $ is
regular and $f^{s}\left(  X^{s}\right)  =f\left(  X\right)  ^{s}.$

\item[(iv)] If $f$ is, in addition, linear and $Z$ is a regular Riesz subspace
of $X$ which is invariant under $f,$ then $Z^{s}$ is invariant under $f^{s}.$
\end{enumerate}
\end{theorem}

\begin{proof}
The proof of (i) is exactly \cite[Proposition 3]{L-444} and (ii) can be proved
in a similar way. Although only existence is proved in \cite[Proposition
3]{L-444}, the ionicity is, however, obvious.

We need only to prove (iii). First we show that $f\left(  X\right)  $ is
regular in $X.$ Consider a net $\left(  x_{a}\right)  _{\alpha\in A}$ in
$f\left(  X\right)  $ such that $x_{\alpha}\uparrow x$ in $f\left(  X\right)
$ and let $y=\sup\limits_{\alpha\in A}x_{\alpha}$ in $X.$ We have to show that
$x=y,$ or equivalently, $x\in f\left(  X\right)  .$ But as $f$ is order
continuous we have $x_{\alpha}=f\left(  x_{\alpha}\right)  \uparrow f\left(
y\right)  $ in $X.$ This shows that $x=f\left(  y\right)  \in f\left(
X\right)  $ as required. Next we will show the equality $f\left(  X\right)
^{s}=f^{s}\left(  X^{s}\right)  .$ We know that if $a\in X^{s},$ then
$f^{s}\left(  a\right)  =\sup f\left(  \left[  a\right]  ^{\leq}\right)  \in
f\left(  X\right)  ^{s}$ (see the proof of \cite[Theorem 6]{L-444}).
Conversely if $u\in f\left(  X\right)  ^{s}$ then again by \cite[Theorem
6]{L-444}, $u=\sup f\left(  A\right)  $ for some upward directed subset $A$ of
$X.$ As $f^{s}$ is left order continuous we get $u=f^{s}\left(  \sup A\right)
\in f^{s}\left(  X^{s}\right)  .$

(iv) This is almost done in the second part of (iii).
\end{proof}

Let $T$ be a conditional expectation operator defined on a Dedekind complete
Riesz space $X$ fixing an order weak unit $e.$ It follows from Theorem
\ref{P1}.(iii) that $R\left(  T\right)  $ a regular Riesz subspace of $X$ and
that $R\left(  T\right)  ^{s}=R\left(  T^{s}\right)  $ where $T^{s}$ denotes
the extension of $T$ to $X^{s}.$ It is easy seen that $T^{s}x=x$ for every
$x\in R\left(  T\right)  ^{s}.$ So we have the following.

\begin{corollary}
Let $T$ be a conditional expectation operator defined on a Dedekind complete
Riesz space $X$ fixing an order weak unit $e.$ Then $R\left(  T\right)  $ is a
regular Riesz subspace of $X$ and $R\left(  T^{s}\right)  =R\left(  T\right)
^{s}.$
\end{corollary}

According to Theorem \ref{P1}, an order projection $P$ defined on $X$ extends
to $X^{s},$ and its extension, denoted by $P^{s},$ satisfies $P^{s}%
x=\sup\limits_{\alpha}Px_{a}$ for every net $\left(  x_{\alpha}\right)  $ in
$X^{s}$ such that $x_{\alpha}\uparrow x.$ The following result supplies Lemma
3.1 and Theorem 3.2 in \cite{L-24}.

\begin{proposition}
With same notations as above, if $g\in\mathcal{R(}T\mathcal{)}_{+}^{s}$, then
\ $\mathcal{R(}T\mathcal{)}^{s}$ is invariant under $P_{g}^{s}$ and $\left(
P_{g}^{d}\right)  ^{s}$ and $B_{g}^{s}$ and $\left(  B_{g}^{d}\right)  ^{s}$
are invariant under $T^{s}.$ Moreover, \ $T^{s}P_{g}^{s}=P_{g}^{s}T^{s}$ \ and
\ $\left(  P_{g}^{d}\right)  ^{s}T^{s}=T^{s}\left(  P_{g}^{d}\right)  ^{s}$.
\end{proposition}

\begin{proof}
We have $P_{g}=P_{P_{g}e},$ and $P_{g}e=\sup ng\wedge e\in R\left(  T\right)
.$ It follows from \cite[Lemme 3.1]{L-24} that $P_{g}T=TP_{g}.$ Observe now
that $P_{g}^{s}T^{s}$ and $T^{s}P_{g}^{s}$ are tow extensions of $P_{g}T$ that
are increasing and left order continuous. The fact that $B_{g}^{s}$ and
$\left(  B_{g}^{d}\right)  ^{s}$ are invariant under $T^{s}$ follows from
Theorem \ref{P1}.(iv). Also the first part of this Theorem shows that
$P_{g}^{s}T^{s}=T^{s}P_{g}^{s}$ as both of them are increasing and left order
continuous and extend $TP_{g}.$ Similarly we show that $\left(  P_{g}%
^{d}\right)  ^{s}T^{s}=T^{s}\left(  P_{g}^{d}\right)  ^{s}$.
\end{proof}

Our next goal is to extend some properties of $X^{s}$ obtained earlier by
Donner or by the first author. We will extend property (P10) by showing that
the Riesz decomposition property is valid for all elements of $X^{s}.$ Recall
that for every $x\in X^{s}$ one can define its positive and negative parts
repectively by $x^{+}=x\vee0$ and $x^{-}=-\left(  x\wedge0\right)  $ and that
$x^{-}\in X.$

\begin{lemma}
Let $x,y,z\in X^{s}$ such that $x\leq y+z.$ \textit{then there exist}
$y_{1},z_{1}\in X^{s}$ \textit{such that}%
\[
y_{1}\leq y,\ z_{1}\leq z\mathit{\ \qquad}\text{\textit{and,}}\mathit{\ \qquad
}x=y_{1}+z_{1}.
\]
If, in addition, $x\in X$ then $y_{1},z_{1}\in X.$
\end{lemma}

\begin{proof}
Assume that $x\leq y+z.$ Then $u=x^{+}+y^{-}+z^{-}\leq y^{+}+z^{+}+x^{-}.$ By
the Riesz decomposition property in $X_{+}^{s},$ (P10), one can write
$u=a+b+c$ with%
\[
0\leq a\leq y^{+},0\leq b\leq z^{+},0\leq c\leq x^{-}.
\]
So $x=y_{1}+z_{1}$ is a required decomposition as $y_{1}=a-y^{-}\leq y$ and
$z_{1}=b-z^{-}+c-x^{-}\leq z.$ Now, if $x\in X$ then $u\in X$ and so $a,b,c\in
X$ by (P4). It follows that $y_{1},z_{1}\in X,$ which completes the proof.
\end{proof}

\begin{corollary}
\label{S6}Let $A$ and $B$ be two subsets of $X^{s}.$ Then $\sup\left(
A+B\right)  =\sup A+\sup B.$
\end{corollary}

\begin{proof}
The inequality $\sup\left(  A+B\right)  \leq\sup A+\sup B$ is obvious. To show
the converse inequality it is enough to prove the following inclusion%
\[
\left[  \sup A+\sup B\right]  ^{\leq}\subseteq\left[  \sup\left(  A+B\right)
\right]  ^{\leq}.
\]
To this end, let $u\in\left[  \sup A+\sup B\right]  ^{\leq}$ and use the Riesz
decomposition property to write $u=x+y$ with $x\leq\sup A$ and $y\in\sup B.$
Pick $a\in A,b\in B$ and observe that%
\[
x\wedge a+y\wedge b\leq a+b\leq\sup\left(  A+B\right)  .
\]
Using properties (P8) and (P6) to get%
\[
\sup\limits_{a\in A,b\in B}\left(  x\wedge a+y\wedge b\right)  =\sup
\limits_{a\in A}\left(  x\wedge a\right)  +\sup\limits_{b\in B}\left(  y\wedge
b\right)  =x+y=u.
\]
We deduce that $u\in\left[  \sup\left(  A+B\right)  \right]  ^{\leq},$ and
this completes the proof.
\end{proof}

\begin{remark}
\label{S17}If $E$ and $F$ are two Dedekind complete Riesz spaces then $E\times
F$ is a Dedekind complete Riesz space with $\left(  E\times F\right)
_{+}=E_{+}\times F_{+}.$ Moreover it is not hard to see that $\left(  E\times
F\right)  ^{s}=E^{s}\times F^{s}$ and then $\left(  E\times F\right)  _{+}%
^{s}=E_{+}^{s}\times F_{+}^{s}.$ This fact combined with Proposition \ref{P3}
will help us to get quick proofs of several properties of $X^{s}$.
\end{remark}

\begin{lemma}
\label{S1}Let $\left(  x_{\alpha}\right)  $ and $\left(  y_{\beta}\right)  $
be two nets in $X_{+}^{s}$ and let $x,y\in X_{+}^{s}.$ such that $x=\sup
x_{\alpha}$ and $y=\sup y_{\beta}.$ Then the following hold.

\begin{enumerate}
\item[(i)] $x+y=\sup\left(  x_{\alpha}+y_{\beta}\right)  .$

\item[(ii)] $x\vee y=\sup\left(  x_{\alpha}\vee y_{\beta}\right)  .$

\item[(iii)] $x\wedge y=\sup\left(  x_{\alpha}\wedge y_{\beta}\right)  .$

\item[(iv)] If $x\wedge y=0$ then $x+y=x\vee y.$

\item[(v)] $\ x+y=x\vee y+x\wedge y.$

\item[(vi)] For every $a\in X_{+}^{s}$ we have
\[
a+x\wedge y=\left(  a+x\right)  \wedge\left(  a+y\right)  ,\text{ and }a+x\vee
y=\left(  a+x\right)  \vee\left(  a+y\right)  .
\]

\end{enumerate}
\end{lemma}

\begin{proof}
(i) follows from (P8) if the nets $\left(  x_{\alpha}\right)  $ and $\left(
y_{\beta}\right)  $ are chosen in $X.$ For the general case we use Corollary
\ref{S6}.

(ii) It is clear that $x\vee y\geq\sup\limits_{\alpha,\beta}x_{\alpha}\vee
y_{\beta}$. On the other hand the inequalities $\sup\limits_{\alpha,\beta
}\left(  x_{\alpha}\vee y_{\beta}\right)  \geq\sup\limits_{\alpha}x_{\alpha
}=x$ and $\sup\limits_{\alpha,\beta}\left(  x_{\alpha}\vee y_{\beta}\right)
\geq y$ are also obvious and then $\sup\limits_{\alpha,\beta}\left(
x_{\alpha}\vee y_{\beta}\right)  \geq x\vee y,$ which gives (ii).

(iii) The inequality $\sup\left(  x_{\alpha}\wedge y_{\beta}\right)  \leq
x\wedge y$ is clear. For the converse we will assume first that $\left(
x_{\alpha}\right)  $ and $\left(  y_{\beta}\right)  $ are in $X.$ Let $u\geq
x_{\alpha}\wedge y_{\beta}$ and let $z\in\left[  x\wedge y\right]  ^{\leq},$
which means that $z\in X,z\leq x$ and $z\leq y.$ Then $u\geq z\wedge
x_{\alpha}\wedge y_{\beta}.$ As this happens for every $\beta$ we get by (P6),
$u\geq z\wedge x_{\alpha}\wedge y=z\wedge x_{\alpha}.$ Using once more (P6) we
obtain $u\geq z\wedge x=z.$ As this happens for every $z\in\left[  x\wedge
y\right]  ^{\leq}$ it follows from (P2) that $u\geq x\wedge y$ and then
$\sup\left(  x_{\alpha}\wedge y_{\beta}\right)  =x\wedge y$.

The general case can be derived from the above case and Theorem \ref{P1}
applied to the map%
\[
f:X_{+}\times X_{+}\longrightarrow X_{+};\left(  x,y\right)  \longmapsto
x\wedge y.
\]
(see Remark \ref{S17}).

(iv) Clearly, \ $x\vee y\leq x+y.$ Let \ $z\in\left[  x+y\right]  ^{\leq}$
then there exist $x^{\prime}\in\left[  x\right]  ^{\leq}$ \ and \ $y^{\prime
}\in\left[  y\right]  ^{\leq}$ \ such that \ $z=x^{\prime}+y^{\prime}.$ But as
$x^{\prime}\bot$ $y^{\prime}$ we get $z=x^{\prime}+y^{\prime}=x^{\prime}\vee
y^{\prime}\leq x\vee y,$ which gives the converse inequality: $x+y\leq x\vee
y.$

(v) and (vi) follow from (i), (ii) and (iii) above and (P2).
\end{proof}

\begin{lemma}
\label{S2}Let $x,y,z\in X_{+}^{s}.$ Then

\begin{enumerate}
\item[(i)] $x\wedge\left(  y+z\right)  \leq x\wedge y+x\wedge z.$

\item[(ii)] If, in addition, $y\wedge z=0$ then $x\wedge\left(  y+z\right)
=x\wedge y+x\wedge z.$

\item[(iii)] If $x\wedge z=0$ then $x\wedge\left(  y+z\right)  =x\wedge y.$
\end{enumerate}
\end{lemma}

\begin{proof}
(i) It follows easily from the Riesz decomposition property in $X^{s}$ that
for every $0\leq u\in\left[  x\wedge\left(  y+z\right)  \right]  ^{\leq}$
there exist $a,b\in X$ such that $0\leq a\leq y,$ $0\leq b\leq z$ and $u=a+b.$
Moreover as $a,b\leq x$ we get $a\leq x\wedge y$ and $b\leq x\wedge z,$ and
then $u\leq x\wedge y+x\wedge z.$ We deduce now from (P2) that
\[
x\wedge\left(  y+z\right)  \leq x\wedge y+x\wedge z.
\]

(ii) In the case when $y\wedge z=0$ we use Lemma \ref{S1}.(iv) to get%
\[
x\wedge(y+z)=x\wedge(y\vee z)=\left(  x\wedge y\right)  \vee\left(  x\wedge
z\right)  =\left(  x\wedge y\right)  +\left(  x\wedge z\right)  .
\]

(iii) This follows easily from the Riesz decomposition property in $X^{s}$ and
property (i) above.
\end{proof}

\begin{corollary}
Let $x,y,v,w\in X_{+}^{s}$ \ such that \ $x\perp y$ \ and \ $v\perp w.$ Then%
\[
(x+v)\wedge(y+w)=x\wedge w\overset{\perp}{+}y\wedge v.
\]

\end{corollary}

\begin{proof}
It is clear that $x\wedge w+y\wedge v$ is a lower bound of $\left\{
x+v,y+w\right\}  .$ On the other hand by Lemma \ref{S2} we have%
\[
\left(  x+v\right)  \wedge\left(  y+w\right)  \leq x\wedge\left(  y+w\right)
+v\wedge\left(  y+w\right)  =x\wedge w+v\wedge y,
\]
which yields the desired equality.
\end{proof}

It is worth noting that if we define a map $\varphi:X\times X\longrightarrow
X$ by putting $\varphi\left(  x\right)  =x\vee y$ then with notation of
Theorem \ref{P1}, the second assertion in Lemma \ref{S1} means that
$\varphi^{s}\left(  x,y\right)  =x\vee y.$ In a similar way properties (i) and
(iii) in the same lemma can be interpreted. The following result gives us
another way to prove the other properties of Lemma \ref{S1}. It will be also
used in the sequel to get quick proofs.

\begin{proposition}
\label{P3}All spaces in the following statements are assumed to be
\textit{Dedekind complete and maps order continuous and increasing.}

\begin{enumerate}
\item[(i)] If $f,g:X\longrightarrow Y$ satisfy $f\leq g$ then $f^{s}\leq
g^{s}.$

\item[(ii)] If $f,g,h:X\longrightarrow Y$ satisfy $h=f+g$ then $h^{s}%
=f^{s}+g^{s}.$

\item[(iii)] If $f:X\longrightarrow Y$ and $g:Y\longrightarrow Z$ then
$\left(  g\circ f\right)  ^{s}=g^{s}\circ f^{s}.$
\end{enumerate}
\end{proposition}

\begin{proof}
The proof of (i) is trivial. (ii) is very similar to Lemma \ref{S1}. Let us
prove (iii). Let $x\in X^{s}$ and let $\left(  x_{\alpha}\right)  $ a net in
$X$ such that $x_{\alpha}\uparrow x.$ Then $f\left(  x_{\alpha}\right)
=f^{s}\left(  x_{\alpha}\right)  \uparrow f^{s}\left(  x\right)  .$ So
$g^{s}\left(  f^{s}\left(  x_{\alpha}\right)  \right)  =g\left(  f\left(
x_{\alpha}\right)  \right)  \uparrow g^{s}\left(  f^{s}\left(  x\right)
\right)  .$ On the other hand $g^{s}\left(  f^{s}\left(  x_{\alpha}\right)
\right)  =\left(  g\circ f\right)  \left(  x_{\alpha}\right)  \uparrow\left(
g\circ f\right)  ^{s}\left(  x\right)  .$ This yields the desired equality
$\left(  g\circ f\right)  ^{s}\left(  x\right)  =g^{s}\left(  f^{s}\left(
x\right)  \right)  .$
\end{proof}

\subsection{Finite and infinite parts of an element in $X^{s}.$\label{FI}}

It was shown in \cite{L-444} that if $Y$ is a regular Riesz subspace of a
Dedekind complete Riesz space $X$ then the subset of $X^{s}$ defined by%
\[
Z=\left\{  \sup A:\emptyset\neq A\subseteq Y\right\}
\]
is the sup-completion of $Y.$ In particular, if $Y$ is an ideal of $X$ then
$Y^{s}\subseteq X^{s}.$ We will denote by $\infty_{Y}$ the greatest element of
$Y.$ If $X=\mathbb{R}^{\Omega}$ is the Riesz space of all real valued
functions defined on a set $\Omega,$ then $X^{s}=\mathbb{R}_{\infty}^{\Omega}$
is consisting of all functions defined on $\Omega$ taking values in
$\mathbb{R}_{\infty}=\mathbb{R}\cup\left\{  \infty\right\}  .$ Every element
in $X^{s}$ can be decomposed as follows: $f=f.\chi_{A}+\infty.\chi
_{\mathbb{R}\backslash A},$ where $A=\left\{  x\in\Omega:f\left(  x\right)
\in\mathbb{R}\right\}  .$ We can say that $f.\chi_{A}$ is the finite part of
$f$ and $\infty.\chi_{\mathbb{R}\backslash A}=f\chi_{\mathbb{R}\backslash A}$
is its infinite part. A similar decomposition can be obtained in the case when
$X=L_{p}\left(  \mu\right)  $ with $1\leq p\leq\infty.$ The goal of this
subsection is to extend this decomposition to a general Dedekind complete
Riesz space.

\begin{theorem}
\label{infinite}Let $x\in X_{+}^{s}.$ Then there is a unique projection band
$B$ in $X$ such that $x=\infty_{B}+u$ with $u\perp B$ in $X^{u}.$
\end{theorem}

\begin{proof}
Let $B$ denote the band of $X$ generated by the subset%
\begin{equation}
\left\{  u\in X:x\geq t\left\vert u\right\vert \text{ for all }t\in\left(
0,\infty\right)  \right\}  , \label{F1}%
\end{equation}
and let $P$ be the corresponding band projection. We extend $P$ and
$P^{d}=I-P$ to $X^{s}.$ It is easily seen that $P^{s}x=\sup H=\sup B$ is the
largest element of $B^{s}$ and that $\left(  P^{d}\right)  ^{s}x\in X.$ Also,
it is not difficult to see that the set defined in (\ref{F1}) is in fact a
band in $X$ and then it coincides with $B.$ Let us denote this set by $H.$ We
will only check that $H$ is closed under addition. To this end, take $a,b$ in
$H.$ Then, by definition, $x\geq2t\left\vert a\right\vert $ and $x\geq
2t\left\vert b\right\vert $ for every $t\geq0.$ It follows that $x\geq
\dfrac{1}{2}\left(  2t\left\vert a\right\vert +2t\left\vert b\right\vert
\right)  \geq t\left\vert a+b\right\vert $ and hence $a+b\in H$ as claimed. We
now observe that for each $a\in B_{+},$ we have%
\[
x\geq P^{s}x\geq P^{s}a=Pa=a.
\]
This shows that $P^{s}x=\infty_{B}.$ We continue to note the band projection
on the band generated by $B$ in $X^{u}$ by $P$ and we claim that $P^{d}x\in
X^{u}.$ Otherwise it follows from \cite[Corollary 15]{L-444} that there exists
$a\in B_{+}^{d}$ such that $P^{d}x\geq ta>0$ for all real $t\geq0.$ This
clearly implies that $a\in H=B,$ a contradiction.

It remains to show that the above decomposition is the unique one. Assume that
$x=\infty_{C}+v$ is another decomposition with $v\perp C^{d}.$ Then for each
$a\in C,$ we have $x=\infty_{C}+v\geq t\left\vert a\right\vert $ for every
$t\in\mathbb{R}_{+}.$ This shows that $C$ is contained in $B$. Conversely if
$a\in B$ then $\left(  P_{C}^{d}\right)  ^{s}x=v\geq tP_{C}^{d}a$ for ever
$t\in\left(  0,\infty\right)  $ and so $P_{C}^{d}a=0,$ as $X$ is Archimedean.
This shows that $a\in C$ and completes the proof.
\end{proof}

\begin{definition}
Let $x\in X_{+}^{s}.$ Then $\infty_{B}$ and $u$ defined in Theorem
\ref{infinite} are the infinite part and finite part of $x$. They will be
denoted by $x^{\infty}$ and $x^{f},$ respectively.
\end{definition}

The following is an immediate consequence of the definition above and it will
be useful later on.

\begin{corollary}
\label{S15}Let $X$ be a Dedekind complete Riesz space and $x\in X_{+}^{s}.$
Then $P\perp P_{x^{\infty}}$ if and only if $Px\in X^{u}.$
\end{corollary}

In \cite{L-444} the first author has shown the following result which provides
a characterization of elements of the sup-completion of a Dedekind complete
Riesz space that do not belong to its universal completion.

\begin{theorem}
\label{T4}(\cite[Theorem 14]{L-444}) Let $X$ be a Dedekind complete Riesz
space with weak order unit $e$ and let $x\in X.$ Then%
\[
x\in X^{s}\setminus X^{u}\Longleftrightarrow\bigwedge\limits_{k=1}^{\infty
}P_{\left(  x-ke\right)  ^{+}}e>0.
\]

\end{theorem}

Now, having defined the finite and infinite parts in $X^{s}$ we can prove a
more precise result.

\begin{theorem}
\label{T5}Let $X$ be a Dedekind complete Riesz space with weak order unit $e$
and let $x\in X_{+}^{s}.$ Then $\bigwedge\limits_{k=1}^{\infty}P_{\left(
x-ke\right)  ^{+}}e=P_{x^{\infty}}e.$
\end{theorem}

\begin{proof}
Write $x=x^{\infty}+x^{f}$ and let $P=P_{x^{\infty}}$ be the band projection
on $X^{u}.$ Observe that%
\[
x-ke=x^{\infty}-kPe+x^{f}-kP^{d}e=x^{\infty}+x^{f}-kP^{d}e.
\]
Thus $\left(  x-ke\right)  ^{+}=x^{\infty}+\left(  x^{f}-kP^{d}e\right)
^{+}.$ Since $x^{\infty}\wedge\left(  x^{f}-kP^{d}e\right)  ^{+}=0,$ it
follows that%
\[
P_{\left(  x-ke\right)  ^{+}}=P_{x^{\infty}}+P_{\left(  x^{f}-kP^{d}e\right)
^{+}}=P_{x^{\infty}}\vee P_{\left(  x^{f}-kP^{d}e\right)  ^{+}},
\]
which implies that%
\[
\bigwedge\limits_{k=1}^{\infty}P_{\left(  x-ke\right)  ^{+}}=P_{x^{\infty}%
}\vee\bigwedge\limits_{k=1}^{\infty}P_{\left(  x^{f}-kP^{d}e\right)  ^{+}%
}=P_{x^{\infty}},
\]
where the last equality follows from Theorem \ref{T4} above as $P^{d}e$ is a
weak order unit in the band $B_{x^{\infty}}^{d}.$
\end{proof}

A more general statement, which can be deduced from Theorem \ref{T5}, is the following:

\begin{corollary}
Let $X$ be a Dedekind complete Riesz space with order weak unit. If $x\in
X_{+}^{s}$ and $u\in X_{+}.$ then
\end{corollary}

\[
\bigwedge\limits_{k=1}^{\infty}P_{\left(  x-ku\right)  ^{+}}=P_{P_{u}%
^{s}x^{\infty}}+P_{\left(  P_{u}^{d}\right)  ^{s}x}=P_{u}P_{x^{\infty}}%
+P_{u}^{d}P_{x}.
\]

The following result is easy to prove and will be used in the proof of Theorem
\ref{T3}.

\begin{proposition}
\label{S16}Let $X$ be a Dedekind complete Riesz space $x,y\in X_{+}^{s}$ and
$\lambda\in\mathbb{R}_{+}$ Then the following hold:

\begin{enumerate}
\item $\left(  x+y\right)  ^{\infty}=x^{\infty}+y^{\infty}$ and $\left(
\lambda x\right)  ^{\infty}=\lambda x^{\infty};$

\item If $x\leq y$ then $x^{\infty}\leq y^{\infty};$

\item $\left(  x\vee y\right)  ^{\infty}=x^{\infty}\vee y^{\infty}$ and
$\left(  x\wedge y\right)  ^{\infty}=x^{\infty}\wedge y^{\infty}.$
\end{enumerate}
\end{proposition}

We end this section with a result which will be needed later in the proof of
Theorem \ref{BCL2}.

\begin{lemma}
\label{S3}Let $X$ be a Dedekind complete Riesz space and let $\left(
x_{n}\right)  $ be a sequence in $X_{+}^{u}.$ Assume that $\sum\limits_{n=1}%
^{\infty}x_{n}=\infty_{B}+u$ as in Proposition \ref{infinite} and let
$R_{n}=\sum\limits_{k=n}^{\infty}x_{k}.$ Then for every $y\in X_{+},$
\[
Py=\lim\left(  y\wedge R_{n}\right)  ,
\]
where $P$ denotes the band projection on $B.$
\end{lemma}

\begin{proof}
Observe first that $P^{d}\sum\limits_{n=1}^{\infty}x_{n}=\sum\limits_{n=1}%
^{\infty}P^{d}x_{n}=u.$ It follows that $P^{d}R_{n}=\sum\limits_{k=n}^{\infty
}P^{d}x_{n}\overset{o}{\longrightarrow}0.$ On the other hand we deduce from
Lemma \ref{S2} that%
\[
y\wedge R_{n}=y\wedge\left(  PR_{n}+P^{d}R_{n}\right)  =y\wedge PR_{n}+y\wedge
P^{d}R_{n}.
\]
As $PR_{n}=\infty_{B}$ for every $n$ and $y\wedge P^{d}R_{n}\overset
{o}{\longrightarrow}0,$ the result follows.
\end{proof}

\subsection{Multiplication in X$^{s}$}

Consider again a Dedekind complete\ Riesz space $X$ with weak order unit $e.$
We know that the universal completion $X^{u}$ of $X$ is equipped with an
$f$-algebra multiplication with $e$ as identity. Our aim in this subsection is
to extend the multiplication to the positive part of the cone $X^{s}$ and
prove several properties of that multiplication which extend standard ones.
Recall, by the way, that although $X^{u}$ is not contained in $X^{s}$ in
general, its positive cone $X_{+}^{u}$ does (\cite[Corollary 7]{L-444}). We
defined the product of two elements $x$ and $y$ of $X_{+}^{s}$ as follows:%
\begin{align*}
xy  &  =\sup\left\{  vw:0\leq v\in\left[  x\right]  ^{\leq},0\leq w\in\left[
y\right]  ^{\leq}\right\} \\
&  =\sup\left\{  vw:0\leq\left(  v,w\right)  \in\left[  \left(  x,y\right)
\right]  ^{\leq}\right\}  .
\end{align*}

Thus the product on $X_{+}^{s}$ is the unique extension of the product on
$X_{+}$ in the sense of Theorem \ref{P1}.(ii). By considering the map
\[
\pi:X_{+}\times X_{+}\longrightarrow X_{+}^{u};\left(  x,y\right)  \longmapsto
xy
\]
one can see that if $x,y\in X_{+}^{+}$ then $xy=\pi^{s}\left(  x,y\right)  .$
The following result is then an immediate consequence of Theorem \ref{P1}.

\begin{lemma}
\label{S5}Let $\left(  x_{\alpha}\right)  ,$ $\left(  y_{\beta}\right)  $ be
two nets in $X_{+}^{s}$ such that $x_{\alpha}\uparrow x$ and $y_{\beta
}\uparrow y.$ Then $x_{\alpha}y_{\beta}\uparrow xy.$
\end{lemma}

As expected the product on $X^{s}$ shares some standard properties.

\begin{lemma}
\label{S7}Let $x,y,z\in X_{+}^{s}.$ The following statements hold.

\begin{enumerate}
\item[(i)] $x\left(  y+z\right)  =xy+xz.$

\item[(ii)] $x\left(  y\wedge z\right)  =xy\wedge xz.$

\item[(iii)] $x\left(  y\vee z\right)  =xy\vee xz.$
\end{enumerate}
\end{lemma}

\begin{proof}
(i) Define the following maps
\begin{align*}
h  &  :X_{+}\times X_{+}\times X_{+}\longrightarrow X_{+};\left(
x,y,z\right)  \longmapsto x\left(  y+z\right)  ;\\
f  &  :X_{+}\times X_{+}\longrightarrow X_{+};\left(  x,y\right)  \longmapsto
xy;\\
g  &  :X_{+}\times X_{+}\longrightarrow X_{+};\left(  x,z\right)  \longmapsto
xz.
\end{align*}
All of them are increasing and order continuous and satisfy $h=f+g.$ Apply
then Proposition \ref{P3} to obtain the equality $h^{s}=f^{s}+g^{s},$ which
means exactly the required equality.

(ii) and (iii) can be proved in a similar way.
\end{proof}

We list next some properties of the product on $X_{+}^{s}.$

\begin{proposition}
Let $X$ be a Dedekind complete Riesz space. $B$ and $C$ be two bands in $X$
and $x,y\in X_{+}^{s}.$ The following properties hold.

\begin{enumerate}
\item[(i)] If $x\in X_{+}^{s}$ and $B\in\mathfrak{B}\left(  X\right)  $ then
$x.\infty_{B}=\infty_{P_{Bx}}.$ In particular, $x.\infty=\infty_{B_{x}}$ and
$\left(  x.\infty_{B}\right)  ^{f}=0.$

\item[(ii)] If $B,C\in\mathcal{B}\left(  X\right)  $ then, $\infty_{B}%
\wedge\infty_{C}=\infty_{B}.\infty_{C}=\infty_{B\cap C}$ and $\infty_{B}%
\vee\infty_{C}=\infty_{B}+\infty_{C}=\infty_{B+C}.$

\item[(iii)] If $x,y\in X_{+}^{s}$ then $\left(  xy\right)  ^{f}=x^{f}y^{f}$
and $\left(  xy\right)  ^{\infty}=x^{\infty}y^{\infty}+x^{\infty}y^{f}%
+x^{f}y^{\infty}.$
\end{enumerate}
\end{proposition}

\begin{proof}
(i) It follows from the definition of the product in $X_{+}^{s}$ that $ex=x$
for all $x\in X_{+}^{s}.$ On other hand Property (P7) yields that if $0<u\in
X$ is a weak order unit then $\sup\limits_{k}ku=\infty.$ Thus%
\[
\infty x=\sup ke.x=\sup\limits_{k\geq1}kx=\infty_{B_{x}}.
\]
for all $x\in X_{+}^{s}.$ So, if $x_{\alpha}\uparrow\infty$ and $x\in
X_{+}^{s}$ then $x_{\alpha}x\uparrow\infty_{B_{x}}.$

(ii) Observe first that for each $n\in\mathbb{N}$ we have%
\[
nP_{B}e\wedge nP_{C}e=n\left(  P_{B}e\wedge P_{C}e\right)  =nP_{B}%
eP_{C}e=nP_{B\cap C}e.
\]
Moreover, as $P_{D}e$ is a weak unit in $D$ for every band $D,$ the result
follows by taking the supremum over $n$ and using Lemma \ref{S1}(iii). For the
second formula observe that%
\[
nP_{B}e\vee nP_{C}e=nP_{B+C}e\leq nP_{B}+nP_{C}e.
\]
By taking the supremum over $n$ we get%
\[
\infty_{B}\vee\infty_{C}=\infty_{B+C}\leq\infty_{B}+\infty_{C}.
\]
On the other hand since $\infty_{B},\infty_{C}\in\left(  B+C\right)  ^{s}$ we
obtain $\infty_{B}+\infty_{C}\leq\infty_{B+C},$ which gives the equality.

(iii) can be deduce easily from (i).
\end{proof}

\section{Borel-Cantelli Lemmas}

We will assume throughout this section that $X$ is a Dedekind complete Riesz
space with conditional expectation operator $T$ and weak order unit $e=Te.$ We
recall that the space $X$ is called $T$-universally complete if every
increasing net $\left(  x_{\alpha}\right)  $ in $X$ with $\left(  Tx_{\alpha
}\right)  $ order bounded in $X^{u}$ is order convergent in $X.$ By extending
$T$ to its natural domain $L^{1}\left(  T\right)  $ we may assume that $X$ is
$T$-universally complete (see \cite{L-24}). The first Borel-Cantelli Lemma has
been generalized to the setting of Riesz spaces in \cite{L-260}. We provide
here a slight more general form of it.

\begin{lemma}
\label{BC1}\textit{Let} $\left(  x_{n}\right)  $ \textit{be an order bounded
sequence of }$X$. \textit{If} $\sum\limits_{n=1}^{\infty}Tx_{n}\in X^{u},$
\textit{then} $\limsup\limits_{n\longrightarrow\infty}x_{n}=0.$
\end{lemma}

\begin{remark}
\begin{enumerate}
If we assume only that $T$ is positive (and not strictly positive), then we
can only conclude that $T\left(  \lim\sup x_{n}\right)  =0.$

\item If $x_{n}=P_{n}e$ for some band projection $P_{n}$ we get a version of
Borel-Cantelli Lemma.\newline If $\sum\limits_{n=1}^{\infty}TPe_{n}\in X^{u},$
\textit{then} $\limsup\limits_{n\longrightarrow\infty}Pe_{n}=0.$
\end{enumerate}
\end{remark}

\begin{proof}
Since $X$ is $T$-universally complete it follows that the series $\sum x_{n}$
is order convergent in $X.$ Now, as $T$ is order continuous we have%
\[
T\sup\limits_{k\geq n}x_{k}\leq T\left(
{\textstyle\sum\limits_{k=n}^{\infty}}
x_{k}\right)  =%
{\textstyle\sum\limits_{k=n}^{\infty}}
Tx_{k},
\]
and then%
\[
T\left(  \limsup\limits_{n\longrightarrow\infty}x_{n}\right)  =\lim
T\sup\limits_{k\geq n}x_{k}=0.
\]
As $T$ is strictly positive we deduce that $\limsup\limits_{n\longrightarrow
\infty}x_{n}=0,$ as required.
\end{proof}

The following lemma is needed in the proof of the second Borel-Cantelli Lemma
in Riesz spaces. Before stating the lemma we recall that the ideal $X_{e}$
generated by $e$ can be endowed with a multiplication for which $e$ is a unit
element. If $x\in X_{e}$ and $f$ is continuous then $f\left(  x\right)  $ is
well defined (\cite[Lemma 2]{L-20}). Moreover it is easily seen that $\left(
fg\right)  \left(  x\right)  =f\left(  x\right)  g\left(  x\right)  .$ One can
also use the $C\left(  K\right)  $-representation of $X_{e},$ when the
constant function $1$ corresponds to the unit $e.$ It should be noted that if
two real functions $f$ and $g$ satisfies $f\leq g$ and $f\left(  x\right)  $
and $g\left(  x\right)  $ exist then $f\left(  x\right)  \leq g\left(
x\right)  .$ As an example of this we get the inequality $e-x\leq\exp\left(
-x\right)  $ for $x\in X_{+},$ which will be needed later. Functional calculus
done on the ideal $X_{e}$ can be extended using Daniell Integral to the whole
of $X.$ This allows to consider more general functions which do not need to be
continuous (see \cite{L-06}). For more information about functional calculus
the reader is referred to \cite{L-06} and \cite{L-180}, where he can find a
comparison between the two kinds mentioned above. In the sequel we need only
the following fact: If $x\in X$ then $\exp\left(  x\right)  $ and $\exp\left(
-x\right)  $ are well defined and $\exp\left(  x\right)  \exp\left(
-x\right)  =e.$

\begin{lemma}
\label{S4}Let $\left(  x_{\alpha}\right)  $ be a net in $X_{+}$ such that
$x_{\alpha}\uparrow x\in X_{+}^{s}.$ Then $\exp\left(  -x_{\alpha}\right)
\downarrow\exp\left(  -x^{f}\right)  $ where $x^{f}$ is the finite part of $x$
defined above. In particular if $x_{\alpha}\uparrow\infty,$ then $\exp\left(
-x_{\alpha}\right)  \downarrow0.$
\end{lemma}

\begin{proof}
Assume first that $x_{\alpha}\uparrow\infty.$ Then $\left(  \exp\left(
-x_{\alpha}\right)  \right)  $ is a positive decreasing net and we have to
show that its infimum is zero. If not, there exists $u>0$ such that
$\exp\left(  -x_{\alpha}\right)  \geq u$ for every $\alpha.$ In this case we
get%
\[
e=\exp\left(  x_{\alpha}\right)  \exp\left(  -x_{\alpha}\right)  \geq
\exp\left(  x_{\alpha}\right)  .u\geq x_{\alpha}u.
\]
Using Lemma \ref{S5} we deduce that $e\geq u^{\infty}.$ But $u^{\infty}\notin
X^{u}$ and this contradicts Property (P4).
\end{proof}

\begin{theorem}
[Borel-Cantelli Lemma]\label{BCL2}Let $\left(  P_{n}\right)  $ be a sequence
of $T$-independent band projections in $X.$ \textit{If} $B$ \textit{is a band
of} $X^{u}$ \textit{such that}%
\[
\sum_{n\geq1}TP_{n}e=\infty_{B}+u\text{ with }u\in B^{d},
\]
\textit{then }$P_{B}$ commutes with $T$ and%
\[
P_{B}=\limsup\limits_{n\longrightarrow\infty}P_{n}.
\]
In particular, if $\sum\limits_{n\geq1}TP_{n}e=\infty$ then $\limsup
\limits_{n\longrightarrow\infty}P_{n}=I.$
\end{theorem}

\begin{proof}
We will show that $Pe\in\mathcal{R}\left(  T\right)  $ and conclude by
\cite[Lemma 3.1]{L-24} that $TP=PT.$ To this end observe first that%
\[
e\wedge\sum\limits_{k=n}^{\infty}TP_{k}e=\lim\limits_{N\longrightarrow\infty
}e\wedge\sum\limits_{k=n}^{N}TP_{k}e\in\mathcal{R}\left(  T\right)  .
\]
According to Lemma \ref{S3} we deduce that $Pe\in\mathcal{R}\left(  T\right)
$ as required.

By the first Borel-Cantelli Lemma (\ref{BC1}) we have $P^{d}\lim\sup
P_{n}e=0.$ So by considering the band $B$ instead of $X$ and $PP_{n}$ instead
of $P_{n}$ we may assume that $B=X$ and $\sum\limits_{n=1}^{\infty}%
TP_{n}e=\infty.$ We have to show that $\limsup P_{n}=I,$ or equivalently
$\left(  \limsup P_{n}\right)  ^{d}e=0.$ Now observe that%
\begin{align*}
\left(  \limsup P_{n}\right)  ^{d}  &  =\liminf P_{n}^{d}=\sup\limits_{n}%
\left(  \inf\limits_{k\geq n}P_{k}^{d}\right) \\
&  =\sup\limits_{n}\left(  \inf\limits_{m\geq n}P_{n}^{d}...P_{m}^{d}\right)
.
\end{align*}
Since the projections $P_{n}$ are $T$-independent we have%
\[
TP_{n}^{d}...P_{m}^{d}e=\prod\limits_{k=n}^{m}\left(  e-TP_{k}e\right)
\leq\exp\left(  -\sum\limits_{k=n}^{m}TP_{k}e\right)  .
\]
Using Lemma \ref{S4} we can see that $TP_{n}^{d}...P_{m}^{d}e\downarrow0$ as
$m\longrightarrow\infty.$ Since all those considering sequences are increasing
or decreasing the order continuity of $T$ allows us to deduce that%
\[
T\left(  \limsup P_{n}\right)  ^{d}e=0.
\]
As $T$ is strictly positive we obtain $\left(  \limsup P_{n}\right)  ^{d}e=0$
as required.
\end{proof}

\section{Applications to Borel Cantelli lemmas}

We will assume also in this section that our space $X$ is Dedekind complete
Riesz space with conditional expectation operator $T$ and weak order unit
$e=Te.$

The following observation may be useful. It could be compared with \cite[Lemma
1.6]{L-421} and \cite[Lemma 1.2]{L-171}.

\begin{lemma}
\label{S12}Assume that $\left(  x_{\alpha}\right)  $ is a positive decreasing
net in $X$ such that $x_{\alpha}\overset{TP}{\longrightarrow}0.$ Then
$x_{\alpha}\downarrow0.$
\end{lemma}

\begin{proof}
Let $x=\inf x_{\alpha}.$ By Lemma \ref{LA}, we have $T\left(  x_{\alpha}\wedge
e\right)  \longrightarrow0.$ On the other hand by order continuity of $T,$
$T\left(  x_{\alpha}\wedge e\right)  \downarrow T\left(  x\wedge e\right)  .$
It follows from the strict positivity of $T$ that $x\wedge e=0$ and then $x=0$
as required.
\end{proof}

Next we present some applications to Borel-Cantelli Lemmas in Riesz spaces. We
need the following lemma.

\begin{lemma}
Let $\left(  x_{n}\right)  _{n\geq1}$ be a sequence in $X$ such that $\left(
\sup\limits_{k\geq0}\left\vert x_{n+k}-x_{n}\right\vert \right)  _{n\geq1}$
converges in $T$-conditional probability to $0$ as $n\longrightarrow\infty,$
then $\left(  x_{n}\right)  _{n\geq1}$ is order convergent in $X^{u}.$
\end{lemma}

\begin{proof}
We have to show that $\delta_{n}=\sup\limits_{p,q\geq n}\left\vert x_{p}%
-x_{q}\right\vert \overset{o}{\longrightarrow}0.$ As order convergence and
uo-convergence agree for sequences in $X^{u}$ (see \cite[Theorem 3.2]{L-63}),
it is sufficient to prove that $\sup\limits_{p,q\geq n}\left\vert x_{p}%
-x_{q}\right\vert \wedge e\overset{o}{\longrightarrow}0$ (\cite[Theoem
28]{L-444}) Now it follows from the assumption and Lemma \ref{LA} that
$T\left(  \sup\limits_{p,q\geq n}\left\vert x_{p}-x_{q}\right\vert \wedge
e\right)  \overset{o}{\longrightarrow}0$ and it follows from Lemma \ref{S12}
that $\sup\limits_{p,q\geq n}\left\vert x_{p}-x_{q}\right\vert \wedge
e\overset{o}{\longrightarrow}0,$ which completes the proof.
\end{proof}

Borel-Cantelli lemma can be used to prove almost surely convergence of
sequences of random variables by showing the convergence of some series.
Similar results can be obtained in the setting of Riesz spaces.

\begin{proposition}
\label{P2}Let $\left(  x_{n}\right)  $ be a sequence in $X$ and let $x\in X.$

\begin{enumerate}
\item[(i)] If the series $%
{\textstyle\sum\limits_{n=1}^{\infty}}
TP_{\left(  \left\vert x_{n}-x\right\vert -\varepsilon e\right)  ^{+}}e$ is
order convergent in $X^{u}$ for all $\varepsilon>0$ then $x_{n}\overset
{uo}{\longrightarrow}x.$

\item[(ii)] If $%
{\displaystyle\sum\limits_{n=1}^{\infty}}
T|x_{n}-x|^{r}\in X^{u}\ $for some real $r>0\ $then $x_{n}\overset
{uo}{\longrightarrow}x.$
\end{enumerate}
\end{proposition}

\begin{proof}
(i) By the first Borel--Cantelli lemma (Lemma \ref{BC1}), we have%
\[
\limsup\limits_{n\longrightarrow\infty}P_{(|x_{n}-x|-\varepsilon e)^{+}%
}e=0\ \ \ \ \varepsilon>0.
\]
So \ $P_{(|x_{n}-x|-\varepsilon e)^{+}}e\overset{o}{\longrightarrow}0.$ The
result now follows from \cite[ Theorem 2.8]{L-63}.

(ii) According to Chebychev Lemma in Riesz space (see \cite[Theorem
3.9]{L-180}), we have
\[
TP_{\left(  \left\vert x_{n}-x\right\vert -\varepsilon e\right)  ^{+}}%
e\leq\varepsilon^{-r}T|x_{n}-x|^{r}.
\]
The result follows now from (i).
\end{proof}

It was proved in \cite{L-750}, that if $\left(  P_{n}\right)  $ is a sequence
of band projections satisfying $P_{n}\uparrow I$ and if $\left(  x_{n}\right)
$ is a sequence in $X$ such that for every $k,$ $\left(  P_{k}x_{n}\right)
_{n\geq1}$ is order convergent in $X^{u},$ then the sequence $\left(
x_{n}\right)  $ is order convergent in $X^{u}.$ We need here a more general statement.

It was shown in \cite{L-750} that if $X$ is a Dedekind complete vector
lattice, $\left(  P_{\gamma}\right)  _{\gamma\in\Gamma}$ is a net of band
projections such that $P_{\gamma}\uparrow I$ and $\left(  x_{\alpha}\right)
_{\alpha\in A}$ is a net of elements in $X$ such that for each $\gamma$ the
net $\left(  P_{\gamma}x_{\alpha}\right)  _{\alpha\in A}$ is order convergent,
then $\left(  x_{\alpha}\right)  $ is $uo$-convergent in $X^{u}.$ Next we
state a more general result: first, the order convergence of the net $\left(
P_{\gamma}x_{\alpha}\right)  _{\alpha\in A}$ is relaxed to $uo$-convergence,
second, the net $\left(  P_{\gamma}\right)  $ is not assumed to be increasing.

\begin{theorem}
\label{S10}Let $\left(  P_{\gamma}\right)  _{\gamma\in\Gamma}$ be a net of
band projections and $\left(  x_{\alpha}\right)  _{\alpha\in A}$ a net of
elements in $X.$ If for each $\gamma,$ the net $\left(  P_{\gamma}x_{\alpha
}\right)  $ is $uo$-convergent in $X,$ then $\left(  Px_{\alpha}\right)  $ is
$uo$-convergent in $X^{u},$ where $P=\sup P_{\gamma}.$
\end{theorem}

\begin{proof}
By considering the band generated by $R\left(  P_{\gamma}\right)  ,$
$\gamma\in\Gamma,$ we may assume that $P=I.$ Let $\mathcal{F}$ be the set of
nonempty finite subsets of $\Gamma.$ For every $F\in\mathcal{F}$ we consider
the band projection $P_{F}=\sup\limits_{\gamma\in F}P_{\gamma}.$ Then $\left(
P_{F}\right)  _{F\in\mathcal{F}}$ is an increasing net of band projections and
$P_{F}\uparrow I.$

Step 1. We will show first that $\left(  P_{F}x_{\alpha}\right)  _{\alpha\in
A}$ is uo-convergent. This can be done by induction on $k=\left\vert
F\right\vert .$ It is exactly what says the assumption of $k=1.$ Now observe
that if $F=G\cup\left\{  \gamma\right\}  $ with $\gamma\notin G$ then%
\[
P_{F}=P_{G}+P_{G}^{d}P_{\gamma}.
\]
As $\left(  P_{\gamma}x_{\alpha}\right)  _{\alpha\in A}$ is uo-convergent, it
follows easily that $\left(  P_{G}^{d}P_{\gamma}x_{\alpha}\right)  _{\alpha\in
A}$ is uo-convergent. So, if we assume that $\left(  P_{G}x_{\alpha}\right)
_{\alpha\in A}$ is $uo$-convergent, then the above equality can be used to
conclude that $\left(  P_{F}x_{\alpha}\right)  _{\alpha\in A}$ is
$uo$-convergent as well. This proves our claim.

Step 2. We will show now that $\left(  x_{\alpha}\right)  _{\alpha\in A}$ is
$uo$-convergent in $X^{u}.$ According to \cite[Theorem 17]{L-444} it is
sufficient to show that $\left(  x_{\alpha}\right)  _{\alpha\in A}$ is
$uo$-Cauchy in $X^{u}$. As $X$ is an ideal in $X^{u}$ and so it is a regular
Riesz subspace, it is even enough to show that $\left(  x_{\alpha}\right)
_{\alpha\in A}$ is $uo$-Cauchy in $X$ (see \cite[Theorem 3.2]{L-65}). By the
first step we know that the net $\left(  P_{F}x_{\alpha}\right)  _{\alpha\in
A}$ is uo-Cauchy, which means that $\left\vert P_{F}x_{\alpha}-P_{F}x_{\beta
}\right\vert \wedge z\overset{o}{\longrightarrow}0$ for each $z\in X_{+}.$ Now
for $y$ in $X_{+}$ we have%
\begin{align*}
\left\vert x_{\alpha}-x_{\beta}\right\vert \wedge y  &  =P_{F}\left(
\left\vert x_{\alpha}-x_{\beta}\right\vert \wedge y\right)  +P_{F}^{d}\left(
\left\vert x_{\alpha}-x_{\beta}\right\vert \wedge y\right) \\
&  =\left\vert P_{F}x_{\alpha}-P_{F}x_{\beta}\right\vert \wedge P_{F}%
y+P_{\gamma}^{d}\left(  \left\vert x_{\alpha}-x_{\beta}\right\vert \wedge
y\right) \\
&  \leq\left\vert P_{F}x_{\alpha}-P_{F}x_{\beta}\right\vert \wedge
P_{F}y+P_{F}^{d}y.
\end{align*}
It follows that%
\[
\limsup\limits_{\left(  \alpha,\beta\right)  }\left\vert x_{\alpha}-x_{\beta
}\right\vert \wedge y\leq P_{F}^{d}y.
\]
As this happens for every $F\in\mathcal{F}$ and $P_{F}^{d}\downarrow0,$ we
deduce that $\limsup\limits_{\left(  \alpha,\beta\right)  }\left\vert
x_{\alpha}-x_{\beta}\right\vert \wedge y=0,$ which proves that $\left(
x_{\alpha}\right)  $ is $uo$-Cauchy in $X$ as claimed. This completes the proof.
\end{proof}

To continue our discussion we need to recall some facts and fix some
notations. Recall that a filtration in $X$ is a sequence of conditional
expectations $\left(  T_{n}\right)  _{n\geq1}$ such that $T_{n}T_{m}%
=T_{m}T_{n}=T_{n\wedge m}$ for every $n,m\geq1.$ We will assume also that
$T_{n}T=TT_{n}=T$ for all $n.$ An adapted process is a sequence $\left(
x_{n}\right)  $ in $X$ such that $x_{n}\in R\left(  T_{n}\right)  $ for every
$n.$ An adapted process is called a martingale if $T_{i}x_{j}=x_{i}$ for all
$i,\ j\in\{1,2,\ \ldots\}$ with $i\leq j$. It is called a submartingale if
$T_{i}x_{j}\geq x_{i}$ for all $i,\ j\in\{1,2,\ \ldots\}$ with $i\leq j.$ We
call a \textit{stopping time\ }adapted to the filtration $(T_{i})_{i\geq1}$ an
increasing sequence $(P_{i})_{i\geq1}$ of band projections on $X$ such that
$P_{i}T_{j}=T_{j}P_{i}$ whenever $1\leq i\leq j$. An example of stopping time
can be obtained as above: if $(x_{k})$ is an increasing sequence in $X$ with
$x_{k}\in R(T_{k})^{+}$ for $k=1,2,\ \ldots$ then the sequence $(P_{k}%
=P_{x_{k}})_{k\geq1}$ is a stopping time. If $P$ is a stopping time, then
$Q=P\wedge n$ defined by
\[
Q_{i}=P_{i}\text{ if }i<n\text{ and }Q_{i}=I\text{ if }i\geq n,
\]
is a stopping time. If $\left(  x_{i}\right)  $ is an adapted process and
$P=\left(  P_{i}\right)  $ is a bounded stopping time, we define the
\textit{stopped process} $(x_{P},T_{P})$ by putting%

\[
x_{P}=\sum_{i=1}^{\infty}\left(  P_{i}-P_{i-1}\right)  x_{i},
\]
where $P_{0}=0.$ We refer the reader to \cite{L-32, L-311} for more
information about the subject.

In the next Lemma and Theorem we consider a submartingale $(x_{n}%
,T_{n})_{n\geq1}.$ For each positive real $K$ we consider a sequence of
projection bands%
\[
B_{K,n}=B\left(  x_{1}\leq Ke,...,x_{n-1}\leq Ke,x_{n}>Ke\right)  .
\]
We use here the notations adopted by Grobler in his papers (see for example
\cite{L-06}), which are close to the notations used in probability theory. So,
$B\left(  x<y\right)  $ denotes the band generated by $\left(  y-x\right)
^{+}$ and $B\left(  x\geq y\right)  $ is its disjoint complement$.$ Thus
$B_{K,n}$ is the intersection of $n$ bands.

The band $B_{K,\infty}$ is defined as follows:%
\[
B_{K,\infty}=B\left(  x_{n}\leq Ke\text{ for all }n\right)  .
\]
Let $P_{K,n}$ denote the corresponding band projection for $n\in\mathbb{N}%
\cup\left\{  \infty\right\}  .$ This allows us to define a stopping time
$\tau^{K}=\left(  \tau_{n}^{K}\right)  _{n\geq1}$ by putting $\tau_{n}^{K}=%
{\textstyle\sum\limits_{j=1}^{n}}
P_{K,j}.$ If $\left(  x_{n}\right)  _{n\geq1}$ is an adapted process then the
stopped process $\left(  x_{n\wedge\tau^{K}}\right)  $ is the process $\left(
z_{n}\right)  _{n\geq1}$ given by
\[
z_{n}=%
{\displaystyle\sum_{j=1}^{n-1}}
P_{_{K,j}}x_{j}+P_{K,n-1}^{d}x_{n}=%
{\displaystyle\sum_{1\leq j\leq\infty}}
P_{_{K,j}}x_{j\wedge n}.
\]

\begin{lemma}
\label{S11}Let $(x_{n},T_{n})_{n\geq1}$ be a submartingale and $K>0.$ If
$\ T[\sup_{n}(x_{n+1}-x_{n})^{+}]$ $\in X^{u}$, \ then the stopped process
$(\widetilde{x}_{n}=x_{n\wedge\tau^{K}})$ satisfies
\[
\limsup\limits_{n\longrightarrow\infty}T|\widetilde{x}_{n}|\in X^{u}.
\]

\end{lemma}

\begin{proof}
We will use the notations preceding Lemma \ref{S11} and recalling that%
\[
\widetilde{x}_{n}=%
{\textstyle\sum\limits_{j=1}^{\infty}}
P_{_{K,j}}x_{j\wedge n}+P_{K,\infty}x_{n}.
\]
Now define a new process $\left(  y_{n}=x_{n\wedge\tau^{K}-1}\right)  $ by
putting%
\[
y_{n}=%
{\textstyle\sum\limits_{j=1}^{\infty}}
P_{_{K,j}}x_{n\wedge j-1}+P_{K,\infty}x_{n},
\]
with $x_{0}=0$ and write%
\[
\widetilde{x}_{n}=y_{n}+\widetilde{x}_{n}-y_{n}.
\]
We have on one hand,%
\[
\widetilde{x}_{n}-y_{n}\leq\left(  \widetilde{x}_{n}-y_{n}\right)  ^{+}\leq
V:=\sup\limits_{j}(x_{j+1}-x_{j})^{+}.
\]
On the other hand as $P_{B_{K,j}}y_{n}^{+}\leq KP_{B_{K,j}}e$ for all
$j\in\mathbb{N\cup}\left\{  \infty\right\}  $ it follows from the definition
of $y_{n}$ that
\[
y_{n}^{+}\leq Ke.
\]
Moreover, as $\left(  \widetilde{x}_{n}\right)  $\ is a submartingale (\cite[
Theorem 4.5]{L-32})\textbf{\ }we have%
\[
T\widetilde{x}_{n}\geq T\widetilde{x}_{1}=Tx_{1}.
\]
As a conclusion, using the above inequalities, we get%
\begin{align*}
T\left\vert \widetilde{x}_{n}\right\vert  &  =2T\widetilde{x}_{n}%
^{+}-T\widetilde{x}_{n}\leq2T\widetilde{x}_{n}^{+}-Tx_{1}\\
&  \leq2Ke+2T\sup_{n}(x_{n+1}-x_{n})^{+}-Tx_{1}\in X^{u},
\end{align*}
which proves the lemma.
\end{proof}

The following result generalizes \cite[Theorem 5.2.8]{b-2435}.

\begin{theorem}
\label{T2}Let $\left(  T_{n}\right)  $ be a filtration adapted with $T$ and
let $\left(  x_{n},T_{n}\right)  $ be a martingale in $X$ satisfying
$T\sup\left\vert x_{n+1}-x_{n}\right\vert \in X^{u}.$ Define the projection
band in $X^{u}$ as follows%
\[
Q=\sup\left\{  P\in\mathfrak{P}:Px_{n}\text{ is order convergent in }%
X^{u}\right\}  .
\]
Then $Q^{d}=P_{\left(  \sup x_{n}\right)  ^{\infty}}.$
\end{theorem}

\begin{proof}
Let $x=\sup\limits_{n}x_{n}.$ By Theorem\textbf{ \ref{S10} }we know that
$\left(  Qx_{n}\right)  $ is uo-convergent in\ $X^{u}.$ Let \ $K\geq1$ \ be a
fixed integer and consider the projection bands $B_{K,n},$ $n\in\mathbb{N}%
\cup\left\{  \infty\right\}  $ and denote by $P_{K,n}$ the corresponding band
projection and let $\tau^{K}$ be the stopping time defined above. It follows
from Lemma \ref{S11} that \ $\limsup\limits_{n\longrightarrow\infty
}T\left\vert \widetilde{x}_{n}\right\vert \in X^{u}.$\textbf{ }This implies by
\cite[Theorem 3.5]{L-03} that\ $(\widetilde{x}_{n})$\ is uo-convergent
in\ $X^{u}.$ Thus $\left(  P_{K,\infty}\widetilde{x_{n}}\right)  $ is
uo-convergent in $X^{u}$ as well. But\ $P_{K,\infty}\widetilde{x}%
_{n}=P_{K,\infty}x_{n}$ and then $P_{K,\infty}x_{n}$ is uo-convergent in
$X^{u}.$ It follows from Theorem \ref{S10} that\textbf{ }$\left(
Px_{n}\right)  $\ is uo-convergent in $X^{u}$ where $P=\sup\limits_{K}%
P_{K,\infty}.$ To prove the theorem it is enough to observe that%
\[
\bigvee\limits_{k=1}^{\infty}B_{K,\infty}=\bigvee\limits_{k=1}^{\infty
}B_{\left(  x-Ke\right)  ^{+}}^{d}=B_{\left(  \sup x_{n}\right)  ^{\infty}%
}^{d}.
\]
Or, equivalently,%
\[%
{\textstyle\bigwedge\limits_{k=1}^{n}}
B_{\left(  x-ke\right)  ^{+}}=B_{x^{\infty}}.
\]
But this is exactly what Theorem \ref{T5} says.
\end{proof}

As a corollary to this theorem we can observe that in the statement of Theorem
5.28 in \cite{b-2435}, $\mathbb{P}\left(  \lim\inf x_{n}=-\infty\right)  =0.$
Furthermore the set $A_{2}$ can be defined as $A_{2}=\left\{  \sup
X_{n}=\infty\right\}  $ instead of $\left\{  \limsup X_{n}=\infty\right\}  .$

The following generalizes L\'{e}vy's Theorem to the frame of Riesz spaces (see
\cite[Theorem 5.1.2]{b-2945}).

\begin{theorem}
\label{T3}Let $\left(  T_{n}\right)  $ be a filtration adapted with $T$ and
let $\left(  P_{n}\right)  $ be a sequence of band projections satisfying
$P_{n}T_{n}=T_{n}P_{n}$ for each integer $n.$ Then the series $%
{\textstyle\sum\limits_{n=1}^{\infty}}
P_{n}e$ and $%
{\textstyle\sum\limits_{n=1}^{\infty}}
T_{n-1}P_{n}e$ have the same infinite part.
\end{theorem}

\begin{proof}
Define for each $n,$ $d_{n}=P_{n}e-T_{n-1}P_{n}e$ and $x_{n}=d_{1}+...+d_{n}.$
Then $\left(  x_{n}\right)  _{n\geq1}$ is a martingale and $\left\vert
x_{n+1}-x_{n}\right\vert =\left\vert d_{n+1}\right\vert \leq e.$ Then Theorem
\ref{T2} can be applied. It is sufficient to show that $Q_{i}%
{\textstyle\sum\limits_{n=1}^{\infty}}
P_{n}e$ and $Q_{i}%
{\textstyle\sum\limits_{n=1}^{\infty}}
T_{n-1}P_{n}e$ have the same infinite part for $i=1,2,$ where $Q_{1}=Q$ and
$Q_{2}=I-Q.$ But this is trivial for $Q_{1}$ because the sequence%
\[
Q_{1}x_{n}=Q_{1}\left[
{\textstyle\sum\limits_{k=1}^{n}}
P_{k}e-%
{\textstyle\sum\limits_{k=1}^{n}}
T_{k-1}P_{k}e\right]
\]
is order convergent in $X^{u}$, and it is also true for $Q_{2}$ as an
immediate consequence of the following inequalities:%
\[
x_{n}\leq%
{\textstyle\sum\limits_{k=1}^{n}}
P_{k}e\qquad\text{and}\qquad-x_{n}\leq%
{\textstyle\sum\limits_{k=1}^{n}}
T_{k-1}P_{k}e,\qquad n\geq1.
\]
The result now follows from Proposition \ref{S16}.
\end{proof}

\end{document}